\documentclass[12pt,reqno]{amsart}
\usepackage{amssymb}
\input amssym.def
\usepackage{amsmath,amsfonts,hyperref,xcolor,mathtools,booktabs}
\usepackage{amscd}
\usepackage[mathscr]{eucal}
\usepackage{enumitem}
\usepackage{orcidlink}
\allowdisplaybreaks
\numberwithin{equation}{section}

\hypersetup{colorlinks=true,citecolor={purple},linkcolor={teal},urlcolor={violet}}

\theoremstyle{plain}
\newtheorem{theorem}{Theorem}[section]

\newtheorem{lemma}[theorem]{Lemma}

\theoremstyle{definition}

\makeatletter
\@namedef{subjclassname@2020}{%
	\textup{2020} Mathematics Subject Classification}
\makeatother

\begin{document}
\title[Remarks on a certain restricted partition function of Lin]{Remarks on a certain restricted partition function of Lin}
\author[Russelle Guadalupe]{Russelle Guadalupe\orcidlink{0009-0001-8974-4502}}
\address{Institute of Mathematics, University of the Philippines Diliman\\
Quezon City 1101, Philippines}
\email{rguadalupe@math.upd.edu.ph}

\renewcommand{\thefootnote}{}



\renewcommand{\thefootnote}{\arabic{footnote}}

\setcounter{footnote}{0}

\begin{abstract}
Let $b(n)$ be the number of partition triples $\pi=(\pi_1,\pi_2,\pi_3)$ of $n$ such that $\pi_1$ consists of distinct odd parts, and $\pi_2$ and $\pi_3$ consist of parts divisible by $4$. Utilizing modular forms, Lin obtained the generating functions for $b(3n+1)$ and $b(3n+2)$, which yields the congruence $b(3n+2)\equiv 0\pmod{3}$ for all $n\geq 0$. We provide in this note elementary proofs of these generating functions by employing $q$-series manipulations and dissection formulas. We also establish infinite families of internal congruences modulo $3$ for $b(n)$.
\end{abstract}

\maketitle

\section{Introduction}\label{sec1}

We denote $f_m := \prod_{n\geq 1}(1-q^{mn})$ for $m\in\mathbb{N}$ and $q\in\mathbb{C}$ with $|q| < 1$ throughout this paper. A partition of a natural number $n$ is a finite nonincreasing sequence of natural numbers, called as its parts, whose sum is $n$. Let $p(n)$ be the number of partitions of $n$ with $p(0):=1$. Then we have the generating function for $p(n)$ given by 
\begin{align*}
\sum_{n\geq 0} p(n)q^n = \dfrac{1}{f_1}.
\end{align*}
Ramanujan \cite[pp. 210--213]{ramp} proved that 
\begin{align}\label{eq11}
\sum_{n\geq 0} p(5n+4)q^n = 5\dfrac{f_5^5}{f_1^6},
\end{align}
which implies that 
\begin{align*}
p(5n+4) \equiv 0\pmod{5}
\end{align*}
for all $n\geq 0$. In connection with the Ramanujan's cubic continued fraction, Chan \cite{chan} defined a cubic partition of $n$ as a partition of $n$ whose even parts may appear in two colors. Letting $a(n)$ be the number of cubic partitions of $n$ with $a(0):=1$, we obtain the generating function
\begin{align*}
\sum_{n\geq 0} a(n)q^n = \dfrac{1}{f_1f_2}.
\end{align*}
Chan \cite{chan} showed that
\begin{align}\label{eq12}
\sum_{n\geq 0} a(3n+2)q^n = 3\dfrac{f_3^3f_6^3}{f_1^4f_2^4},
\end{align}
which is an analogue of (\ref{eq11}) and follows that
\begin{align*}
a(3n+2)\equiv 0\pmod{3}
\end{align*}
for all $n\geq 0$. Subsequently, Kim \cite{kim} introduced an overcubic partition of $n$ as a cubic partition of $n$ in which the first occurrence of each part may be overlined. Letting $\overline{a}(n)$ be the number of overcubic partitions of $n$ with $a(0):=1$, we get the generating function
\begin{align*}
\sum_{n\geq 0} \overline{a}(n)q^n = \dfrac{f_4}{f_1^2f_2}.
\end{align*}
With the help of modular forms, Kim \cite{kim} found the following identity
\begin{align}\label{eq13}
\sum_{n\geq 0} \overline{a}(3n+2)q^n = 6\dfrac{f_3^6f_4^3}{f_1^8f_2^3}
\end{align}
analogous to (\ref{eq11}) and (\ref{eq12}), yielding
\begin{align*}
\overline{a}(3n+2)\equiv 0\pmod{6}
\end{align*}
for all $n\geq 0$. Hirschhorn \cite{hirsc} provided an elementary proof of (\ref{eq13}) and derived the generating functions for $\overline{a}(3n)$ and $\overline{a}(3n+1)$.

In 2013, Lin \cite{lin} explored arithmetic properties of the restricted partition function $b(n)$ that counts the number of partition triples of $\pi=(\pi_1,\pi_2,\pi_3)$ of $n$ such that $\pi_1$ consists of distinct odd parts, and $\pi_2$ and $\pi_3$ consist of parts divisible by $4$. The generating function for $b(n)$ with $b(0):=1$ is given by 
\begin{align}\label{eq14}
\sum_{n\geq 0} b(n)q^n = \dfrac{f_2^2}{f_1f_4^3}.
\end{align}
Lin \cite{lin} employed modular forms to prove the following generating functions.
\begin{theorem}[\cite{lin}]\label{thm11}
We have
\begin{align}
	\sum_{n\geq 0} b(3n+2)q^n &= 3q\dfrac{f_2^6f_{12}^6}{f_1^3f_4^{11}},\label{eq15}\\
	\sum_{n\geq 0} b(3n+1)q^n &= \alpha(q^4)\dfrac{f_2^6f_{12}^3}{f_1^3f_4^{10}},\label{eq16}
\end{align}
where $\alpha(q)$ is the cubic theta function \cite{bbg}
\begin{align*}
	\alpha(q) := \sum_{m, n=-\infty}^{\infty} q^{m^2+mn+n^2}.
\end{align*}
\end{theorem}

We note that (\ref{eq15}) resembles (\ref{eq12}) and (\ref{eq13}) and implies the congruence 
\begin{align}\label{eq17}
b(3n+2)\equiv 0\pmod{3}
\end{align}
for all $n\geq 0$. We also note that there is a quick proof of (\ref{eq17}), as we now illustrate. Applying the $3$-dissection \cite[(14.3.3)]{hirscb}
\begin{align}\label{eq18}
\dfrac{f_2^2}{f_1} = \dfrac{f_6f_9^2}{f_3f_{18}}+q\dfrac{f_{18}^2}{f_9}
\end{align}
and the congruence $f_k^3\equiv f_{3k}$ for $k\geq 1$ by the binomial theorem (which will be frequently used without further comment) on (\ref{eq14}), we arrive at
\begin{align*}
\sum_{n\geq 0} b(n)q^n \equiv \dfrac{f_2^2}{f_1f_{12}}\equiv \dfrac{1}{f_{12}}\left(\dfrac{f_6f_9^2}{f_3f_{18}}+q\dfrac{f_{18}^2}{f_9}\right)\pmod{3}.
\end{align*}
Comparing the terms involving $q^{3n+2}$ in the above congruence leads to (\ref{eq17}) as desired. 

The objective of this paper is to give an elementary proof of Theorem \ref{thm11} by relying on the classical $q$-series manipulations and dissection formulas. We also explore the arithmetic properties of the generating function for $b(3n)$ and derive the following families of internal congruences modulo $3$ as a consequence.

\begin{theorem}\label{thm12}
For all $n\geq 0$ and $k\geq 0$, we have 
\begin{align}
	b\left(3^{2(k+1)}n+\dfrac{5\cdot 3^{2(k+1)}+3}{8}\right)&\equiv (-1)^k b(9n+6)\pmod{3},\label{eq19}\\
	b(81n+24)&\equiv -b(9n+3)\pmod{3}.\label{eq110}
\end{align}	
\end{theorem}

We lay out the rest of the paper as follows. Section \ref{sec2} is devoted to showing Theorem \ref{thm11} using certain identities involving $\alpha(q)$, $2$-dissections, and $3$-dissections. We prove in Section \ref{sec3} Theorem \ref{thm12} by applying $3$-dissections of the generating functions for certain restricted partition functions due to Toh \cite{toh} and Andrews, Hirschhorn, and Sellers \cite{ahs}.

\section{Another Proof of Theorem \ref{thm11}}\label{sec2}

We first require the following identities necessary to show (\ref{eq15}) and (\ref{eq16}). We have the following $2$-dissections \cite[Lemma 1]{hirscr}
\begin{align}
\dfrac{f_3^3}{f_1} &= \dfrac{f_4^3f_6^2}{f_2^2f_{12}}+q\dfrac{f_{12}^3}{f_4},\label{eq21}\\
\dfrac{f_1}{f_3^3} &= \dfrac{f_2^3f_{12}^3}{f_4f_6^9}\left(\dfrac{f_4^3f_6^2}{f_2^2f_{12}}-q\dfrac{f_{12}^3}{f_4}\right),\label{eq22}
\end{align}
and the following $3$-dissection \cite[Lemma]{hirsca}
\begin{align}\label{eq23}
\dfrac{1}{f_1^3} = \dfrac{f_9^3}{f_3^{10}}\left(\alpha(q^3)^2+3q\alpha(q^3)\dfrac{f_9^3}{f_3}+9q^2\dfrac{f_9^6}{f_3^2}\right).
\end{align}

We also require the following identities \cite[(22.11.6), (22.6.1)]{hirscb}
\begin{align}
\alpha(q) &= \dfrac{f_2^6f_3}{f_1^3f_6^2} + 3q\dfrac{f_6^6f_1}{f_3^3f_2^2},\label{eq24}\\
\alpha(q^4) &= \alpha(q) - 6q\dfrac{f_4^2f_{12}^2}{f_2f_6}.\label{eq25}
\end{align}

We next deduce the following important theta function identities.

\begin{lemma}\label{lem21}
We have 
\begin{align*}
	\dfrac{f_6^8f_1}{f_3^4f_2^2}+\dfrac{f_2f_3^2f_{12}^3}{f_1f_4f_6} = 2\dfrac{f_4^2f_{12}^2f_6}{f_2f_3}.
\end{align*}	
\end{lemma}

\begin{proof}
Using (\ref{eq21}) and (\ref{eq22}), we write
\begin{align*}
	\dfrac{f_6^8}{f_2^2}\cdot\dfrac{f_1}{f_3^3}+\dfrac{f_2f_{12}^3}{f_4f_6}\cdot\dfrac{f_3^3}{f_1}
	&=\dfrac{f_6^8}{f_2^2}\cdot\dfrac{f_2^3f_{12}^3}{f_4f_6^9}\left(\dfrac{f_4^3f_6^2}{f_2^2f_{12}}-q\dfrac{f_{12}^3}{f_4}\right)+\dfrac{f_2f_{12}^3}{f_4f_6}\left(\dfrac{f_4^3f_6^2}{f_2^2f_{12}}+q\dfrac{f_{12}^3}{f_4}\right)\\
	&=\dfrac{f_2f_{12}^3}{f_4f_6}\cdot 2\dfrac{f_4^3f_6^2}{f_2^2f_{12}}=2\dfrac{f_4^2f_6f_{12}^2}{f_2}.
\end{align*}
Dividing both sides by $f_3$ yields the desired identity.
\end{proof}

\begin{lemma}\label{lem22}
We have 
\begin{align*}
	\alpha(q^4)\dfrac{f_6^2}{f_3}+3q\dfrac{f_2f_3^2f_{12}^3}{f_1f_4f_6} = \dfrac{f_2^6}{f_1^3}.
\end{align*}	
\end{lemma}

\begin{proof}
We infer from (\ref{eq24}) and (\ref{eq25}) that
\begin{align*}
	\alpha(q^4)\dfrac{f_6^2}{f_3} +3q\dfrac{f_2f_3^2f_{12}^3}{f_1f_4f_6}
	&=\dfrac{f_6^2}{f_3}\left(\dfrac{f_2^6f_3}{f_1^3f_6^2} + 3q\dfrac{f_6^6f_1}{f_3^3f_2^2}- 6q\dfrac{f_4^2f_{12}^2}{f_2f_6}\right)+3q\dfrac{f_2f_3^2f_{12}^3}{f_1f_4f_6}\\
	&=\dfrac{f_2^6}{f_1^3}+3q\left(\dfrac{f_6^8f_1}{f_3^4f_2^2}-2\dfrac{f_4^2f_{12}^2f_6}{f_2f_3}+\dfrac{f_2f_3^2f_{12}^3}{f_1f_4f_6}\right)\\
	&=\dfrac{f_2^6}{f_1^3},
\end{align*}
where the last equality follows from Lemma \ref{lem21}.
\end{proof}

We are now in a position to give another proof of Theorem \ref{thm11}.

\begin{proof}[Proof of Theorem \ref{thm11}]
We apply (\ref{eq18}) and (\ref{eq23}) on (\ref{eq14}) to obtain
\begin{align}
	\sum_{n\geq 0} b(n)q^n = \left(\dfrac{f_6f_9^2}{f_3f_{18}}+q\dfrac{f_{18}^2}{f_9}\right)\cdot \dfrac{f_{36}^3}{f_{12}^{10}}\left(\alpha(q^{12})^2+3q^4\alpha(q^{12})\dfrac{f_{36}^3}{f_{12}}+9q^8\dfrac{f_{36}^6}{f_{12}^2}\right).\label{eq26}
\end{align}
We then extract the terms of (\ref{eq26}) involving $q^{3n+2}$, divide both sides by $q^2$, and replace $q^3$ with $q$. In view of Lemma \ref{lem22}, we get
\begin{align*}
	\sum_{n\geq 0} b(3n+2)q^n &= 3q\dfrac{f_6^2}{f_3}\cdot \alpha(q^4)\dfrac{f_{12}^6}{f_4^{11}} + 9q^2\dfrac{f_2f_3^2}{f_1f_6}\cdot \dfrac{f_{12}^9}{f_4^{12}}\\
	&= 3q\dfrac{f_{12}^6}{f_4^{11}}\left(\alpha(q^4)\dfrac{f_6^2}{f_3}+3q\dfrac{f_2f_3^2f_{12}^3}{f_1f_4f_6}\right)\\
	&= 3q\dfrac{f_2^6f_{12}^6}{f_1^3f_4^{11}},
\end{align*}
which is precisely (\ref{eq15}). We next consider the terms of (\ref{eq26}) involving $q^{3n+1}$, divide both sides by $q$, and replace $q^3$ with $q$. We see from Lemma \ref{lem22} that
\begin{align*}
	\sum_{n\geq 0} b(3n+1)q^n &= \alpha(q^4)^2\dfrac{f_6^2}{f_3}\cdot \dfrac{f_{12}^3}{f_4^{10}}+3q\alpha(q^4)\dfrac{f_2f_3^2}{f_1f_6}\cdot \dfrac{f_{12}^6}{f_4^{11}}\\
	&=\alpha(q^4)\dfrac{f_{12}^3}{f_4^{10}}\left(\alpha(q^4)\dfrac{f_6^2}{f_3}+3q\dfrac{f_2f_3^2f_{12}^3}{f_1f_4f_6}\right)\\
	&=\alpha(q^4)\dfrac{f_2^6f_{12}^3}{f_1^3f_4^{10}},
\end{align*}
which is exactly (\ref{eq16}).
\end{proof}

\section{Proof of Theorem \ref{thm12}}\label{sec3}

We need certain $3$-dissections to establish Theorem \ref{thm12}. We have the identities \cite[(14.3.2)]{hirscb}, \cite[Lemma 2.1, (2.1c)]{toh}
\begin{align}
\dfrac{f_1^2}{f_2} &= \dfrac{f_9^2}{f_{18}}-2q\dfrac{f_3f_{18}^2}{f_6f_9},\label{eq31}\\
\dfrac{f_2}{f_1f_4} &= \dfrac{f_{18}^9}{f_3^2f_9^3f_{12}^2f_{36}^3}+q\dfrac{f_6^2f_{18}^3}{f_3^3f_{12}^3}+q^2\dfrac{f_6^4f_9^3f_{36}^3}{f_3^4f_{12}^4f_{18}^3}.\label{eq32}
\end{align}
We note that (\ref{eq32}) is the $3$-dissection of the generating function for the number of partitions of $n$ with distinct odd parts, which was found by Toh \cite{toh}. On the other hand, Andrews, Hirschhorn,
and Sellers \cite[Theorem 3.1]{ahs} obtained the $3$-dissection
\begin{align}\label{eq33}
\dfrac{f_4}{f_1} = \dfrac{f_{12}f_{18}^4}{f_3^3f_{36}^2}+q\dfrac{f_6^2f_9^3f_{36}}{f_3^4f_{18}^2}+2q^2\dfrac{f_6f_{18}f_{36}}{f_3^3}
\end{align}
of the generating function for the number of partitions of $n$ with distinct even parts. 

We next present the following theta function identity needed in the proof of Theorem \ref{thm12}.

\begin{lemma}\label{lem31}
We have 
\begin{align*}
	\dfrac{f_2^2f_3f_{12}^3}{f_4^4f_6^2}+\dfrac{f_2^4f_3^5f_{12}^3}{f_1^4f_4^4f_6^4}=2\dfrac{f_2f_6^5}{f_1^2f_3f_4^3}.
\end{align*}
\end{lemma}

\begin{proof}
We require the following $2$-dissection \cite[(22.6.2)]{hirscb}, \cite[(1.35)]{hirscgb}
\begin{align}
	\dfrac{f_1^3}{f_3} &= \dfrac{f_4^3}{f_{12}}-3q\dfrac{f_2^2f_{12}^3}{f_4f_6^2}.\label{eq34}
\end{align}
Using (\ref{eq21}), (\ref{eq22}), and (\ref{eq34}), we compute
\begin{align}
	\dfrac{f_2f_{12}^3}{f_4f_6^2}&\cdot f_1^2f_3^2+\dfrac{f_2^3f_{12}^3}{f_4f_6^4}\cdot\left(\dfrac{f_3^3}{f_1}\right)^2 = \dfrac{f_3^3}{f_1}\left[\dfrac{f_2f_{12}^3}{f_4f_6^2}\cdot \dfrac{f_1^3}{f_3}+\dfrac{f_2^3f_{12}^3}{f_4f_6^4}\cdot\dfrac{f_3^3}{f_1}\right]\nonumber\\
	&=\dfrac{f_3^3}{f_1}\left[\dfrac{f_2f_{12}^3}{f_4f_6^2}\left(\dfrac{f_4^3}{f_{12}}-3q\dfrac{f_2^2f_{12}^3}{f_4f_6^2}\right)+\dfrac{f_2^3f_{12}^3}{f_4f_6^4}\left(\dfrac{f_4^3f_6^2}{f_2^2f_{12}}+q\dfrac{f_{12}^3}{f_4}\right)\right]\nonumber\\
	&=2\dfrac{f_3^3}{f_1}\left(\dfrac{f_2f_4^2f_{12}^2}{f_6^2}-q\dfrac{f_2^3f_{12}^6}{f_4^2f_6^4}\right)\nonumber\\
	&=2\dfrac{f_3^3}{f_1}\cdot f_6^5 \cdot \dfrac{f_2^3f_{12}^3}{f_4f_6^9}\left(\dfrac{f_4^3f_6^2}{f_2^2f_{12}}-q\dfrac{f_{12}^3}{f_4}\right)\nonumber\\
	&=2\dfrac{f_3^3}{f_1} \cdot f_6^5\cdot \dfrac{f_1}{f_3^3} = 2f_6^5. \label{eq35}
\end{align}
Multiplying both sides of (\ref{eq35}) by $f_2/(f_1^2f_3f_4^3)$, we obtain the desired identity.
\end{proof}

We are now ready to prove Theorem \ref{thm12}.

\begin{proof}[Proof of Theorem \ref{thm12}]	
We first deduce from (\ref{eq24}) and (\ref{eq25}) that $\alpha(q^4)\equiv \alpha(q)\equiv 1\pmod{3}$. Thus, extracting the terms of (\ref{eq26}) involving $q^{3n}$, replacing $q^3$ with $q$, and applying (\ref{eq32}) yield
\begin{align}
	\sum_{n\geq 0} b(3n)q^n &\equiv \dfrac{f_2f_3^2}{f_1f_6}\cdot \alpha(q^{4})^2\dfrac{f_{12}^3}{f_4^{10}} \equiv \dfrac{f_3^2}{f_6}\cdot \dfrac{f_2}{f_1}\cdot \dfrac{f_{12}^3}{f_4f_4^9}\nonumber\\
	&\equiv \dfrac{f_3^2}{f_6}\cdot \dfrac{f_2}{f_1f_4}\nonumber\\
	&\equiv \dfrac{f_3^2}{f_6}\left(\dfrac{f_{18}^9}{f_3^2f_9^3f_{12}^2f_{36}^3}+q\dfrac{f_6^2f_{18}^3}{f_3^3f_{12}^3}+q^2\dfrac{f_6^4f_9^3f_{36}^3}{f_3^4f_{12}^4f_{18}^3}\right)\pmod{3}.\label{eq36}
\end{align}
We consider the terms of (\ref{eq36}) involving $q^{3n+2}$, divide both sides by $q^2$, and then replace $q^3$ with $q$. In view of (\ref{eq31}) and (\ref{eq32}), we have that
\begin{align}
	\sum_{n\geq 0} b(9n+6)q^n &\equiv \dfrac{f_1^2}{f_2}\cdot \dfrac{f_2^4f_3^3f_{12}^3}{f_1^4f_4^4f_6^3}\equiv \dfrac{f_1}{f_4}\cdot \dfrac{f_3^2f_{12}^2}{f_6^2}\label{eq37}\\
	&\equiv \dfrac{f_3^2f_{12}^2}{f_6^2}\left( \dfrac{f_9^2}{f_{18}}-2q\dfrac{f_3f_{18}^2}{f_6f_9}\right)\nonumber\\
	&\times \left(\dfrac{f_{18}^9}{f_3^2f_9^3f_{12}^2f_{36}^3}+q\dfrac{f_6^2f_{18}^3}{f_3^3f_{12}^3}+q^2\dfrac{f_6^4f_9^3f_{36}^3}{f_3^4f_{12}^4f_{18}^3}\right)\pmod{3}.\label{eq38}
\end{align}
We look for the terms of (\ref{eq38}) involving $q^{3n+2}$, divide both sides by $q^2$, and then replace $q^3$ with $q$. We infer from Lemma \ref{lem31} and (\ref{eq33}) that
\begin{align}
	\sum_{n\geq 0} b(27n+24)q^n&\equiv \dfrac{f_1^2f_4^2}{f_2^2}\left(-2\dfrac{f_2f_6^5}{f_1^2f_3f_4^3}+\dfrac{f_2^4f_3^5f_{12}^3}{f_1^4f_4^4f_6^4}\right)\nonumber\\
	&\equiv \dfrac{f_1^2f_4^2}{f_2^2}\cdot -\dfrac{f_{12}^3f_2^2f_3}{f_4^4f_6^2}\equiv -\dfrac{f_4}{f_1}\cdot \dfrac{f_3^2f_{12}^2}{f_6^2}\nonumber\\
	&\equiv -\dfrac{f_3^2f_{12}^2}{f_6^2}\left(\dfrac{f_{12}f_{18}^4}{f_3^3f_{36}^2}+q\dfrac{f_6^2f_9^3f_{36}}{f_3^4f_{18}^2}+2q^2\dfrac{f_6f_{18}f_{36}}{f_3^3}\right)\pmod{3}.\label{eq39}
\end{align}
We now consider the terms of (\ref{eq39}) involving $q^{3n+1}$, divide both sides by $q$, and then replace $q^3$ with $q$. We deduce that
\begin{align}\label{eq310}
	\sum_{n\geq 0} b(81n+51)q^n\equiv -\dfrac{f_1^2f_4^2}{f_2^2}\cdot \dfrac{f_2^2f_3^3f_{12}^2}{f_1^4f_6^2}\equiv -\dfrac{f_1}{f_4}\cdot \dfrac{f_3^2f_{12}^2}{f_6^2}.
\end{align}
Comparing (\ref{eq37}) and (\ref{eq310}), we arrive at
\begin{align*}
	b(81n+51)\equiv -b(9n+6)\pmod{3}
\end{align*}
for all $n\geq 0$. Hence, (\ref{eq19}) follows from the above congruence and induction on $k$.

On the other hand, we look for the terms of (\ref{eq36}) involving $q^{3n+1}$, divide both sides by $q$, and then replace $q^3$ with $q$. We get
\begin{align}\label{eq311}
	\sum_{n\geq 0} b(9n+3)q^n\equiv \dfrac{f_1^2}{f_2}\cdot\dfrac{f_2^2f_6^3}{f_1^3f_4^3}\equiv \dfrac{f_2f_6^3}{f_1f_{12}}\pmod{3}.
\end{align}
We now consider for the terms of (\ref{eq39}) involving $q^{3n}$ and replace $q^3$ with $q$, yielding
\begin{align}\label{eq312}
	\sum_{n\geq 0} b(81n+24)q^n\equiv -\dfrac{f_1^2f_4^2}{f_2^2}\cdot \dfrac{f_4f_6^4}{f_1^3f_{12}^2}\equiv -\dfrac{f_2f_6^3}{f_1f_{12}}\pmod{3}.
\end{align}
Comparing (\ref{eq311}) and (\ref{eq312}), we finally obtain (\ref{eq110}) as desired.
\end{proof}

\section*{Declaration of interests}

The author declares that he has no competing interests.

\section*{Author contributions}

The author conceptualized, wrote and reviewed the entire manuscript.

\section*{Data availability}

There is no data involved in this work.

\end{document}